\documentclass[a4paper]{article}
\usepackage[T1]{fontenc}
\usepackage[english]{babel}
\usepackage{csquotes}
\usepackage[a4paper, hcentering, verbose]{geometry}
\usepackage{xcolor}
\definecolor{TUCGreen}{HTML}{005F50} % 228B22; TU Chemnitz: 005F50

\usepackage{amsthm, mathtools, amssymb, nicefrac}
\mathtoolsset{showonlyrefs}
\usepackage[bookmarks=true, bookmarksnumbered=false, pdfborder={0 0 1}, backref=false, pdfauthor={Alois Pichler}, pdftitle={multistage trees}, colorlinks=true, citecolor=TUCGreen, urlcolor=darkgray, linkcolor=TUCGreen]{hyperref}
\usepackage{doi}

\usepackage{dsfont}%	für \one und Erwartungswert, doublestroke
\usepackage{microtype}
\usepackage[ruled, lined, linesnumbered]{algorithm2e} % | noend
\usepackage{orcidlink}% after hyperref

\DeclareMathOperator{\E}{{\mathds E}}

\DeclareMathOperator{\supp}{supp}		% variance

\DeclareMathOperator*{\minimize}{minimize}
\newcommand{\nd}{\mathsf{d\mkern-1.0mu{l}}}

\numberwithin{equation}{section}
\theoremstyle{plain}
	\newtheorem{theorem}             {Theorem}[section]
	\newtheorem{corollary}  [theorem]{Corollary}
	
	\newtheorem{proposition}[theorem]{Proposition}
\theoremstyle{definition}
	\newtheorem{definition} [theorem]{Definition}
\theoremstyle{remark}
	\newtheorem{remark}     [theorem]{Remark}
\usepackage[round]{natbib}	% References: round, numbers
\title{\textbf{Approximation of Multistage Stochastic Processes}\\ In Memory of Roger Jean-Baptiste\ Wets}
\author{Georg Ch.~Pflug\thanks{University of Vienna, \orcidlink{0000-0001-8215-3550} \protect\protect\href{http://orcid.org/0000-0001-8215-3550}{orcid.org/0000-0001-8215-3550}\protect}
	\and Alois Pichler\thanks{University of Technology, Chemnitz, Faculty of Mathematics, 09126~Chemnitz, Germany \\
	\orcidlink{0000-0001-8876-2429} \href{http://orcid.org/0000-0001-8876-2429}{orcid.org/0000-0001-8876-2429}\protect \\
	Contact: \href{mailto:alois.pichler@math.tu-chemnitz.de}{alois.pichler@math.tu-chemnitz.de}}}
\begin{document}
\maketitle
% Journal: Communications in Optimization Theory (COT), https://cot.mathres.org/
\begin{abstract}
	Multistage stochastic optimization problems are driven by stochastic processes that model uncertain future outcomes. Because exact computation is practically intractable, the underlying process is replaced by a finite scenario tree, and the reliability of the resulting optimal decisions depends on how well the tree represents the original process.

	We measure the approximation quality by the nested distance —~the natural multiperiod extension of the Wasserstein distance~— which provides a rigorous bound on the suboptimality incurred when the true process is replaced by its tree approximation. We study two scenario tree constructions. For trees built by stage-wise optimal quantization, we determine the branching numbers that achieve a prescribed accuracy at minimal total scenario count.
	For randomly sampled trees, we derive an exponential concentration inequality for the nested distance via a union bound over a classical empirical-measure concentration result, and optimize the allocation of samples across stages.
	The two constructions are complementary: quantization is more scenario-efficient in low dimensions, random sampling in high dimensions, and each yields the correct scaling of the approximation error with the number of scenarios.
\end{abstract}

\noindent\textbf{Classification:} 62G05, 62G08, 62G20

\section{Introduction and Mathematical Setting}
Multistage stochastic programming was pioneered by \citet{Dantzig1955} and \citet{Wets1966}, who established the complete-recourse framework that underlies the formulations considered here.
Scenario trees form the computational backbone of numerical methods for multistage stochastic programming \citep{RuszczynskiShapiro2009}.
Replacing the true probability measure~$P$ of a stochastic process by a finite discrete tree~$\tilde P$ is unavoidable in practice, yet every such replacement introduces an error that propagates into the optimal decisions.
Quantifying and controlling this error is therefore fundamental to the reliability of stochastic programming \citep{DupachovaWets1988, DupachovaGroweKuschaRomisch2003}.

For two-stage problems, the Wasserstein distance is the natural error metric: by the Kanto\-ro\-vich--Rubinstein theorem \citep{Villani2003} --~going back to Kantorovich in 1942~-- , the error in the optimal objective is bounded by the Wasserstein distance between~$P$ and~$\tilde P$, weighted by the Lipschitz constant of the cost \citep{Romisch2003, DupachovaGroweKuschaRomisch2003}.
Optimal quantization theory \citep{Zador, GrafLuschgy} then determines the minimum number of support points needed to achieve a prescribed accuracy.

The multistage case is fundamentally more demanding.
Non-anticipative policies require that the decision at each stage~$t$ may depend only on the information available up to that stage \citep{Rockafellar1976, WetsRockafellar97}, so approximating the marginals of a process well does not suffice: the conditional distributions must be approximated simultaneously, consistently with the information, i.e. the filtration structure.
The asymptotic stability of optimal solutions under perturbation of the probability measure was established by \citet{DupachovaWets1988}; the nested distance provides a quantitative version of this stability.
The \emph{nested distance}, introduced by \citet{Pflug2009} and further developed in \citet{PflugPichler2011, PflugPichlerBuch}, extends the Wasserstein distance to account for this filtration structure, and plays the same role for multistage problems as the Wasserstein distance plays for two-stage ones.
Scenario trees and the aggregation of policies in stochastic programs were studied systematically by \citet{RockafellarWets1991}.
For background on variational convergence and epi-convergence we refer to \citet{WetsRockafellar97}.

This paper derives quantitative bounds on the nested distance between a stochastic process and its scenario-tree approximations, and determines the optimal branching structure of such trees. 
Two complementary constructions are developed.
The first is a \emph{deterministic quantization tree}: at each node the conditional distribution is replaced by its optimal quantization, yielding a nested-distance guarantee that holds with certainty, and we determine the branching numbers that achieve a prescribed accuracy at minimal total scenario count.
The second is a \emph{random sampling tree}: at each node independent samples are drawn from the conditional kernel, and under a sub-Gaussian concentration condition \citep{Bolley}, the nested distance is small with high probability; optimal allocation of samples across stages then determines the branching numbers for this probabilistic guarantee.
The two constructions are complementary: quantization is more scenario-efficient when the stage dimensions are small, random sampling when they are large.
\bigskip

The remainder of this section recalls the two-stage setting, optimal quantization, the nested distance, and derives the optimal deterministic branching structure (Proposition~\ref{prop:branches}).

Let $\xi= (\xi_1,\dots,\xi_T)$ be a stochastic process in discrete time, $t= 1,\dots,T$ with values in $\Xi$.
The stochastic optimization problem is
\begin{equation}
	\minimize_{x(\cdot)\in\mathcal{X}} \E c\bigl(\xi,x(\xi)\bigr), \label{eq:1}
\end{equation}
where $c$ is a cost function and $\mathcal{X}$ collects the feasible non-anticipative policies, i.e., measurable functions $x=(x_1, \dots, x_T)$ adapted to the filtration $\mathcal{F}$ generated by~$\xi$.
The expectation is taken with respect to the probability measure~$P$ on $(\Xi,\mathcal{F})$.

Evaluating the expectation in~\eqref{eq:1} for varying policies requires integration over a high-dimensional probability space, which is generally intractable.
We therefore replace~$P$ by a discrete measure~$\tilde P$ supported on a finite scenario tree, for which the expectation reduces to a finite sum.
The key question is how the optimal value and optimal policy of the approximate problem relate to those of the original: this requires Lipschitz continuity of the optimal value as a function of the underlying measure, with the nested distance as the metric on the space of measures.

In what follows we recall the Wasserstein distance on a metric space $(\Xi,d)$ and the two-stage setting, then extend both to the multistage case via the nested distance.

Section~\ref{sec:random} introduces random scenario trees, establishes a probabilistic bound on the nested distance via a union-bound argument (Theorem~\ref{thm:random}), and identifies the optimal branching for random trees (Proposition~\ref{prop:erroralloc}).

\subsection{Two-stage stochastic programs}
We first address two-stage stochastic optimization problems.
In this setting, the feasible policies are $m^\prime$-vectors, that is, the problem~\eqref{eq:1} is
\begin{equation}\label{eq:2}
	\minimize_{x\in\mathcal{X}}\E c(\xi,x), 
\end{equation} where $\mathcal{X}\subset \mathbb{R}^{m^\prime}$.
To quantify how well~$\tilde P$ approximates~$P$, we need a metric on the space of probability measures that is sensitive to the geometry of~$\Xi$ and compatible with the expected-cost objective in~\eqref{eq:2}.
The Wasserstein distance provides exactly this: it measures the minimum work required to transport one distribution into the other, and directly controls the error in the optimal value.
\begin{definition}[Wasserstein distance]
	Let $P$ and $\tilde{P}$ be Borel probability measures on the metric space $(\Xi,d)$.
	The Wasserstein distance of~$P$ and~$\tilde{P}$ of order~$r$ ($r\ge1$) with respect to the metric $d\colon \Xi\times\Xi\to[0,\infty)$ is
	\begin{equation}\label{eq:6}
		w_r(P,\tilde{P})\coloneqq \inf_\pi\left(\iint_{\Xi\times\Xi} d(x,y)^r\,\pi(\mathrm{d}x,\mathrm{d}y)\right)^{\nicefrac{1}r},
	\end{equation}
	where the infimum is over all couplings (transport plans, i.e.\ bivariate measures)~$\pi$ with marginals~$P$ and~$\tilde{P}$, that is,
	\begin{align*}
		\pi(A\times \Xi) &= P(A),\quad\text{ and}\\
		\pi(\Xi\times B) &= \tilde{P}(B)
	\end{align*}
	for every set $A$, $B\in\mathcal{B}(\Xi)$, the Borel $\sigma$\nobreakdash-algebra on~$\Xi$.
\end{definition}
We refer to \citet{Villani2003} for a comprehensive treatment of the Wasserstein distance and its properties. The couplings~$\pi$ in~\eqref{eq:6} are transport plans that redistribute mass from~$P$ to~$\tilde P$ at minimum average cost; the infimum is always attained and the resulting~$w_r$ is a metric on the space of probability measures with finite $r$-th moment \citep{Villani2009}.
The central role of the Wasserstein distance in stochastic optimization stems from the following duality result, which bounds the error in any expected cost when the underlying measure is perturbed.

\begin{proposition}[Kantorovich--Rubinstein theorem]
	For any Lipschitz function $g\colon\Xi\to \mathbb R$,
	\begin{equation}
		|\E_P g- \E_{\tilde P} g|\le \operatorname{Lip}(g)\cdot w_{1}(P,\tilde{P}),
	\end{equation}
	where $\operatorname{Lip}(g)$ denotes the Lipschitz constant of~$g$ and $\E_P g= \int_\Xi g(\xi)\,P(\mathrm{d}\xi)$ is the expectation of~$g$ with respect to the measure~$P$ ($\E_{\tilde P}g= \int_\Xi g(\xi)\,\tilde{P}(\mathrm{d}\xi)$, resp.).
\end{proposition}
Applying the theorem to $g(\cdot)= c(\,\cdot\,,x)$ for all decisions $x\in\mathcal{X}$, the error in the optimal value of~\eqref{eq:2} when $P$ is replaced by $\tilde P$ is bounded by $\operatorname{Lip}(c)\cdot w_1(P,\tilde P)$, where $\operatorname{Lip}(c)$ is a Lipschitz constant of the cost in its first argument, uniformly for the second argument. A small Wasserstein distance between~$P$ and~$\tilde P$ therefore guarantees a small error in both the optimal value and the optimal decisions \citep{DupachovaWets1988, Romisch2003}.
The question then becomes: how efficiently can~$\tilde P$ be chosen, and how many support points are required to achieve a prescribed accuracy?
This is the classical problem of \emph{optimal quantization}, to which we now turn.

\subsection{Optimal quantization}
In what follows, we recall a central result from \citet{Zador} and \citet[Theorem~6.2]{GrafLuschgy} (cf.\ also \citet{Na}) on optimal quantization on $\Xi=\mathbb{R}^{m}$. Optimal quantization seeks to approximate the measure~$P$ by a finitely supported measure $\tilde P= \sum_{i=1}^n p_{i}\cdot\delta_{x_{i}}$ with $n$ support points $\{x_{1},\dots,x_n\}\subset \Xi$ and weights $p_{i}\ge 0$ with $\sum_{i=1}^n p_{i}=1$.

By the Lebesgue decomposition theorem, $P= P_a+ P_s$ where~$P_a$ is absolutely continuous and~$P_s$ is singular with respect to the Lebesgue measure; we write $f_P(x)= \frac{P_a(\mathrm{d}x)}{\mathrm{d}x}$ for the density of the absolutely continuous part.
Zador’s theorem then states that the optimal Wasserstein approximation error 
satisfies the asymptotic relation
\begin{equation}\label{eq:Zador}
	n^{\nicefrac1m}\cdot\inf_{\substack{p_{i}\ge 0,\ \sum_{i=1}^n p_{i}=1,\\ x_{i}\in\Xi,\ i=1,\dots,n}} w_r\biggl(P,\ \sum_{i=1}^n p_{i}\cdot\delta_{x_{i}}\biggr)
			\xrightarrow[n\to\infty]  {} Q_{m,r}\cdot\left\|f_P\right\|^{\nicefrac{1}{r}}_{\frac{m}{m+r}}, 
\end{equation}
where $Q_{m,r}$ is the \emph{quantization coefficient}, a geometric constant depending only on the dimension~$m$ and order~$r$, defined as the asymptotic quantization error for the uniform distribution on~$[0,1]^m$.
By the Kantorovich--Rubinstein theorem, this directly controls the optimal-value error: for each fixed $x\in\mathcal{X}$,
\[
	\left|\int_\Xi c(\xi,x)\,P(\mathrm{d}\xi) -\int_\Xi c(\xi,x)\,\tilde P(\mathrm{d}\xi)\right| \le \operatorname{Lip}(c(\cdot,x))\cdot w_1(P,\tilde P) \le \operatorname{Lip}(c(\cdot,x))\cdot w_r(P,\tilde P),
\]
so the worst-case value error is at most
$\sup_{x\in\mathcal{X}}\operatorname{Lip}(c(\cdot,x))\cdot w_r(P,\tilde P)$.

Setting $q_{m,r}\coloneqq Q_{m,r}\cdot\left\|f_P\right\|^{\nicefrac{1}{r}}_{\frac{m}{m+r}}$, the optimal quantization error behaves asymptotically as
\[	w_r\left(P,\,\sum_{i=1}^{n}p_i\,\delta_{x_i}\right)\;\sim\; q_{m,r}\cdot n^{-1/m},\]
and the number of support points required to achieve accuracy~$\varepsilon$ is
\[	n_\varepsilon \sim \left(\frac{q_{m,r}}\varepsilon\right)^{m}.\]
Since $n_\varepsilon$ grows as $\varepsilon^{-m}$, optimal quantization suffers from the \emph{curse of dimensionality}.

The optimal support points can be found by minimizing the non-linear and non-convex mapping
\[	(x_1,\dots,x_n)\mapsto \int_\Xi\min_{i=1,\dots,n} d(x_i,\xi)^r\,P(\mathrm{d}\xi),\]
and the optimal weights are then the Voronoi cell probabilities
\[	p_i\coloneqq P\left(\bigl\{\xi\colon d(x_i,\xi)=\min_{j=1,\dots,n}  d(x_j,\xi)\bigr\}\right).\]

The analysis above covers the two-stage case, in which feasible policies are constants.
For multistage problems, policies must adapt to the observed history at each stage, and an approximation that is close in Wasserstein distance may still produce large errors in optimal decisions if it fails to reproduce the conditional structure of the process.
The appropriate metric for the multistage setting is the \emph{nested distance}, which we recall next.

\subsection{Multistage stochastic programs}
The two-stage stochastic problem~\eqref{eq:2} is a special case with a nonrandom decision set $\mathcal{X}\subset \mathbb{R}^{m^\prime}$.
Multistage stochastic optimization generalizes this by allowing policies of the form 
\begin{equation}\label{eq:3}
	x(\xi)= \begin{pmatrix*}[l]
		~x_{0}\\
		~x_{1}(\xi_{1})\\
		~x_{2}(\xi_{1},\xi_{2})\\
		~\vdots\\
		~x_T(\xi_{1},\xi_{2},\dots,\xi_T)
	\end{pmatrix*}. 
\end{equation}
The stochastic process $\xi=(\xi_1,\dots,\xi_T)$ has components~$\xi_t$ taking values in $\Xi_t\subset\mathbb{R}^{m_t}$.

The decision problem takes the form \eqref{eq:1} with policies as in~\eqref{eq:3} above.
The observations~$\xi_t$ in time~$t$ take values in the metric space $\Xi_t \subset \mathbb{R}^{m_t}$ with metric~$d_t$; we write $\xi_{:t}\coloneqq(\xi_1,\dots,\xi_t)$ for the history up to stage~$t$.
The decisions $x_t$ lie in $\mathcal{X}_t \subseteq \mathbb{R}^{m_t^\prime}$. At stage~$t$, the policy~$x_t \in \mathbb{R}^{m_t^\prime}$ depends only on the observations $(\xi_1,\dots,\xi_t)$ available up to that stage; this form in~\eqref{eq:3} is the mathematical expression of time-consistency. The distance on $\Xi= \Xi_1\times \dots \times \Xi_T$ is 
\begin{equation}
	d(x,y)= \sum_{t=1}^T d_t(x_t,y_t).	\label{eq:dist1}
\end{equation}

The nested distance extends the Wasserstein distance to stochastic processes and plays the same role for the multistage problem~\eqref{eq:1} as the Wasserstein distance does for the two-stage problem~\eqref{eq:2}.
The nested distance was introduced in \citet{Pflug2009}.

\begin{definition}[Nested distance]
	Let $P$ and $\tilde{P}$ be probability measures on a filtered measurable space $\bigl(\Xi,\mathcal{F}= (\mathcal{F}_t)_{t=1}^T\bigr)$, where $\mathcal{F}_t$ is the $\sigma$\nobreakdash-algebra generated by the coordinate projection $\pi_t\colon \Xi\to \Xi_{:t}\coloneqq \Xi_1\times\cdots\times\Xi_t$.
	The \emph{nested distance} with cost function~$d$ is
	\[	\nd_r(P,\tilde{P})\coloneqq \left(\inf_\pi \iint_{\Xi\times\Xi}d(x,y)^{r}\,\pi(\mathrm{d}x,\mathrm{d}y)\right)^{\nicefrac{1}{r}}, \] where the infimum is over all bivariate measures~$\pi$ on $\Xi\times\Xi$ such that, for every $t=1,\dots,T$ and $\pi$\nobreakdash-almost every $(\xi_{:t-1},\tilde\xi_{:t-1})$,
	\begin{align*}
		\pi(A\times \Xi\mid\xi_{:t-1},\tilde\xi_{:t-1}) &= P(A\mid\xi_{:t-1})\text{ and}\\
		\pi(\Xi\times B\mid\xi_{:t-1},\tilde\xi_{:t-1}) &= \tilde{P}(B\mid\tilde\xi_{:t-1}),
	\end{align*}
	for all $A$, $B\in\mathcal{B}(\Xi_t)$, where $\pi(\cdot\mid\cdot)$ is the conditional distribution of $(\xi_t,\tilde\xi_t)$ given the joint history.
\end{definition}

\begin{remark}
The motivation for the name of this distance lies the fact that it can be constructed in a backward recursion, where for each stage the original distance of the state space is replaced by the distance of the corresponding subtrees \citep{Pflug2023}: the distance for later stages is \emph{nested} within the distance for earlier stages.  
\end{remark}

\begin{remark}[Mapping upper bound]\label{rem:1}
	When $P$ admits conditional densities and $\tilde P$ is a discrete measure, the nested distance can be bounded via \emph{consistent transport mappings}.
	A sequence $S_t\colon \Xi_{:t}\to\Xi_{:t}$, $t=1,\dots,T$, is \emph{consistent} if
	\begin{equation}\label{eq:consistent}
		\operatorname{proj}_s\circ S_t = S_s \qquad \text{for all } s< t,
	\end{equation}
	and each $S_t$ pushes $P_{:t}$ forward to $\tilde P_{:t}$.
	For any such family,
	\begin{align}\label{eq:mapping_bound}
		\nd_r(P,\tilde P) \le \left(\int_\Xi d\bigl(\xi,\,S_T(\xi)\bigr)^r P(\mathrm{d}\xi)\right)^{1/r},
	\end{align}
	since the coupling $\pi= (\mathrm{id},S_T)_\# P$ is feasible for the nested distance: consistency~\eqref{eq:consistent} implies that its conditional marginals match $P_t(\,\cdot\mid\xi_{:t-1})$ and $\tilde P_t(\,\cdot\mid\tilde\xi_{:t-1})$ at every stage. The consistency condition is the mapping-level analog of non-anticipativity.
\end{remark}

In what follows, we shall concentrate on the product $\Xi\coloneqq \Xi_{1}\times\dots\times\Xi_T$ with $\Xi_t= \mathbb{R}^{m_t}$, and the filtration generated by the natural projection.
By the disintegration theorem (cf.~\citet{Ambrosio2005}), any measure~$P$ on the product~$\Xi$ may be disintegrated so that
\begin{equation}
	P(A_1\times\dots\times A_T)= \int_{A_1}\int_{A_{2}}\ldots\int_{A_T}P_T(\mathrm{d}\xi_T\mid \xi_1,\dots,\xi_{T-1})\ \ldots\ P_2(\mathrm{d}\xi_2\mid \xi_1)\ P_1(\mathrm{d}\xi_1\mid\xi_0),
\end{equation} or, in compact notation,
\begin{equation}\label{eq:disintegrated}
	P(\mathrm{d}\xi_1,\dots,\mathrm{d}\xi_T) = P_1(\mathrm{d}\xi_1\mid\xi_0)\,P_2(\mathrm{d}\xi_2\mid\xi_1)\,\cdots\,P_T(\mathrm{d}\xi_T\mid\xi_{:T-1}),
\end{equation} where $P_t(\,\cdot\mid \xi_1,\dots,\xi_{t-1})$ is the kernel measure at stage~$t$ and $\xi_0$ is a deterministic starting value for the scenario process $\xi$, typically the today's known value.
We have the following result on approximations for the nested distance (cf.\ \citet{Mirkov2007b}). 

\begin{theorem}[Submartingale-like bound]\label{thm:Mirkov}
	Let $P$ and $\tilde{P}$ be probability measures on~$\Xi$.
	Suppose there are stage-wise functions $D_t\colon \Xi_{:t}\times\Xi_{:t}\to [0,\infty)$, $t=0,\dots,T$, with $D_{0}\ge0$ constant and $D_T=d$, the distance. The functions $D_t$ serve as stage-wise cost functions and need not themselves be metrics.
	For histories $\xi_{:t-1}$, $\tilde{\xi}_{:t-1}$ set 
	\begin{equation}\label{eq:4}
		\mathcal{W}_t(\xi_{:t-1},\tilde{\xi}_{:t-1})^{r}\coloneqq \inf_{\pi_t}\iint_{\Xi_t\times\Xi_t} D_t\bigl((\xi_{:t-1},\xi_t),(\tilde{\xi}_{:t-1},\tilde{\xi}_t)\bigr)^{r}\,\pi_t(\mathrm{d}\xi_t,\mathrm{d}\tilde{\xi}_t),
	\end{equation}
	where $\pi_t$ ranges over all couplings of $P_t(\,\cdot\mid\xi_{:t-1})$ and $\tilde{P}_t(\,\cdot\mid\tilde{\xi}_{:t-1})$.

	Suppose there are deterministic constants $\varepsilon_t$, $\kappa_t\ge0$ such that for $t=1, \dots, T$
	\begin{equation}\label{eq:cond}
		\mathcal{W}_t(\xi_{:t-1},\tilde{\xi}_{:t-1})\le \varepsilon_t+\kappa_t\,D_{t-1}(\xi_{:t-1},\tilde{\xi}_{:t-1})\quad\text{for all }\xi_{:t-1}\text{ and }\tilde{\xi}_{:t-1}
	\end{equation}
	holds uniformly over all histories at each stage $t=1,\dots,T$. Suppose further that the starting values $\xi_0$ and $\tilde{\xi}_0$ are identical.

	Then the nested distance is bounded by
	\begin{equation}\label{eq:5}
		\nd_r(P,\tilde{P})\le \sum_{t=1}^T\varepsilon_t\prod_{j=t+1}^T\kappa_{j}.
	\end{equation}
\end{theorem}
The condition~\eqref{eq:cond} states that the conditional Wasserstein distance $\mathcal{W}_t$ between the kernels grows at most linearly in the history distance $D_{t-1}$, with Lipschitz constant~$\kappa_t$ and uniform bias~$\varepsilon_t$.
\begin{remark}
	Each term $\varepsilon_t\prod_{j=t+1}^T\kappa_j$ captures the approximation error introduced at stage~$t$, amplified by the expansion of the process over the remaining stages. In particular, if the process is non-expanding ($\kappa_t=1$ for all $t$), the bound simplifies to
	\[	\nd_r(P,\tilde{P}) \le \sum_{t=1}^T\varepsilon_t. \]
\end{remark}

\begin{remark}
	We may set $D_t=d$ and for $t < T$ 
	\[	D_t(\xi_{:t},\tilde{\xi}_{:t}) \coloneqq w_r\left(P_{t+1}(\,\cdot\,\mid\xi_{:t}),\; \tilde{P}_{t+1}(\,\cdot\,\mid\tilde{\xi}_{:t});\,D_{t+1}\right), \]
	i.e., the Wasserstein distance between the stage-$(t+1)$ kernels with ground cost~$D_{t+1}$, which is of course a metric. Notice that by this construction $D_{t+1}$ is nested in $D_t$.
	For the optimal coupling~$\pi_{t+1}^*$ of $P_{t+1}(\,\cdot\,\mid\xi_{:t})$ and $\tilde{P}_{t+1}(\,\cdot\,\mid\tilde{\xi}_{:t})$, we then have
	\[	D_t(\xi_{:t},\tilde{\xi}_{:t})^r = \iint_{\Xi_{t+1}\times\Xi_{t+1}} D_{t+1}\bigl((\xi_{:t},\xi_{t+1}),\,(\tilde{\xi}_{:t},\tilde{\xi}_{t+1})\bigr)^r \,\pi_{t+1}^*(\mathrm{d}\xi_{t+1},\mathrm{d}\tilde{\xi}_{t+1}). \]
\end{remark}
\begin{remark}
	Suppose $r=1$ and the stage-wise cost takes the weighted form, slightly more general than~\eqref{eq:dist1},
	\[	D_t(\xi_{:t},\tilde\xi_{:t}) = \sum_{s=1}^{t}\beta_s\cdot d_s(\xi_s,\tilde\xi_s), \quad \beta_s>0, \]
	where $d_s$ is the basic metric  on $\Xi_s$.  Since the terms
	$\beta_1\,d_1(\xi_1,\tilde\xi_1),\dots,\beta_{t-1}\,d_{t-1}(\xi_{t-1},\tilde\xi_{t-1})$
	are constant in the $\pi_t$-integration, the quantity defined in~\eqref{eq:4} satisfies
	\[
		\mathcal{W}_t(\xi_{:t-1},\tilde\xi_{:t-1})
			= D_{t-1}(\xi_{:t-1},\tilde\xi_{:t-1}) + \beta_t\,w_1\bigl(P_t(\cdot\mid\xi_{:t-1}),\, \tilde{P}_t(\cdot\mid\tilde\xi_{:t-1});\,D_t\bigr).
	\]
	Hence, condition~\eqref{eq:cond} at stage~$t$ reduces to
	\[
		w_1\bigl(P_t(\cdot\mid\xi_{:t-1}),\, \tilde{P}_t(\cdot\mid\tilde\xi_{:t-1});\,D_t\bigr)
			\le \frac{\varepsilon_t}{\beta_t} + \frac{\kappa_t-1}{\beta_t} \sum_{i=1}^{t-1}\beta_i\,d_i(\xi_i,\tilde\xi_i).
	\]
	In particular, for $\kappa_t=1$ the bound is the constant $\varepsilon_t/\beta_t$, independent of the history.
\end{remark}

\begin{proof}[Proof of Theorem~\ref{thm:Mirkov}]
	For each $t$, by a measurable selection theorem (cf.~\citet[Chapter~5]{Villani2003}), there exists a measurably chosen optimal coupling $\pi_t(\,\cdot\mid\xi_{:t-1},\tilde\xi_{:t-1})$ in the definition of~$\mathcal{W}_t$.
	Glue them with respect to their history and set (cf.\ the disintegrated measure in~\eqref{eq:disintegrated})
	\[
		\pi(A\times B)\coloneqq \iint_{A_{1}\times B_{1}}\dots\iint_{A_T\times B_T}\pi_T(\mathrm{d}\xi_T,\mathrm{d}\tilde{\xi}_T\mid \xi_{:T-1},\tilde{\xi}_{:T-1})\dots\pi_{1}(\mathrm{d}\xi_{1},\mathrm{d}\tilde{\xi}_{1}),
	\]
	where $A=A_{1}\times\dots\times A_T$ and $B=B_{1}\times\dots\times B_T$.
	By construction, the conditional marginals of~$\pi$ given $(\xi_{:t-1},\tilde{\xi}_{:t-1})$ are $P_t(\,\cdot\mid\xi_{:t-1})$ and $\tilde{P}_t(\,\cdot\mid\tilde{\xi}_{:t-1})$, so~$\pi$ is feasible for the nested distance.
	
	Define $I_t\coloneqq \| \mathcal W_t\|_{L^r(\pi)} = \bigl(\E_\pi\,D_t(\xi_{:t},\tilde{\xi}_{:t})^{r}\bigr)^{\nicefrac1r}$.
	By condition~\eqref{eq:cond} and Minkowski’s inequality,
	\[	I_t= \|\mathcal{W}_t\|_{L^r(\pi)}\le \|\varepsilon_t+ \kappa_t\,D_{t-1}\|_{L^r(\pi)} \le \varepsilon_t+ \kappa_t\,\|D_{t-1}\|_{L^r(\pi)}= \varepsilon_t+ \kappa_t\,I_{t-1}. \]
	Unrolling the recursion with $I_0= 0$ gives
	\[	I_T\le  \sum_{t=1}^T \varepsilon_t \prod_{j=t+1}^T\kappa_j. \]
	Since $\pi$ is feasible for the nested distance and $D_T = d$,
	\[	\nd_r(P,\tilde{P}) \le \Bigl(\iint d(\xi,\tilde\xi)^r\,\pi(\mathrm{d}\xi,\mathrm{d}\tilde\xi)\Bigr)^{\nicefrac{1}{r}} \le  \sum_{t=1}^T\varepsilon_t\prod_{j=t+1}^T\kappa_j, \]
	which is~\eqref{eq:5}.
\end{proof}

\begin{proposition}\label{prop:Delta}
	Suppose the measure~$P$ satisfies
	\[	w_r\bigl(P_{t+1}(\,\cdot\mid\xi_{:t}),\;P_{t+1}(\,\cdot\mid\tilde{\xi}_{:t})\bigr) \le \kappa_t\cdot D_t(\xi_{:t},\tilde{\xi}_{:t}) \]
	with expansion coefficient $\kappa_t$ for all histories $\xi_{:t},\tilde{\xi}_{:t}\in\Xi_{:t}$.
	Then
	\[	w_r\bigl(P_{t+1}(\,\cdot\mid\xi_{:t}),\;\tilde{P}_{t+1}(\,\cdot\mid\tilde{\xi}_{:t})\bigr) \le \varepsilon_{t+1}+\kappa_t\cdot D_t(\xi_{:t},\tilde{\xi}_{:t}), \]
	where $\varepsilon_{t+1}\coloneqq \sup_{\tilde{\xi}_{:t}} w_r\bigl(P_{t+1}(\,\cdot\mid\tilde{\xi}_{:t}),\; \tilde{P}_{t+1}(\,\cdot\mid\tilde{\xi}_{:t})\bigr)$, which verifies condition~\eqref{eq:cond} at stage~$t+1$.
\end{proposition}
\begin{proof}
	By the triangle inequality of the Wasserstein distance, 
	\begin{align*}
		\MoveEqLeft[5] w_r\left(P_{t+1}(\,\cdot\mid\xi_{:t}),\ \tilde{P}_{t+1}(\,\cdot\mid\tilde{\xi}_{:t})\right)\\
	  &\le w_r\left(P_{t+1}(\,\cdot\mid\xi_{:t}),\ P_{t+1}(\,\cdot\mid\tilde{\xi}_{:t})\right)+ w_r\left(P_{t+1}(\,\cdot\mid\tilde{\xi}_{:t}),\ \tilde{P}_{t+1}(\,\cdot\mid\tilde{\xi}_{:t})\right)\\
					& \le\kappa_t\cdot D_t(\xi_{:t},\tilde{\xi}_{:t})+ w_r\left(P_{t+1}(\,\cdot\mid\tilde{\xi}_{:t}),\ \tilde{P}_{t+1}(\,\cdot\mid\tilde{\xi}_{:t})\right)\\
			& \le \varepsilon_{t+1}+ \kappa_t\cdot D_t(\xi_{:t},\tilde{\xi}_{:t}),
	\end{align*}
	the assertion.
\end{proof}

While Theorem~\ref{thm:Mirkov} bounds the nested distance through the stage-wise errors~$\varepsilon_t$, the following proposition provides an alternative estimate expressed directly as an integrated transport cost, involving only transport mapping (cf.\ Remark~\ref{rem:1}).
\begin{proposition}[Stage-wise mapping decomposition]\label{prop:mapping}
	Suppose that $P$ admits conditional densities.
	For each history $\xi_{:t-1}$, let $S_{t}(\xi_t\mid \xi_{:t-1})$ be a  transport map from $P_t(\,\cdot\mid\xi_{:t-1})$ to $\tilde P_t(\,\cdot\mid\xi_{:t-1})$ under~$d_t$.
	Define recursively $\mathcal{S}_1 (\xi_1) = S_1(\xi_1)$ and
	$\mathcal{S}_t(\xi_{:t})= \bigl(S_1(\xi_1), S_2(\xi_2 \mid \xi_1), \dots, S_t(\xi_t \mid \mathcal{S}_{t-1}(\xi_{:t-1}))\bigr)$.  This family is consistent by construction. It then follows that
	\begin{equation}
		\nd_r(P,\tilde P)\le \left(\int_\Xi d\bigl(\xi, \mathcal{S}_T(\xi)\bigr)^r\;P(\mathrm{d}\xi)\right)^{\nicefrac1r}
			\le \sum_{t=1}^T \left(\int_{\Xi_{:t}} d_t\bigl(\xi_{:t}, \mathcal{S}_t(\xi_{:t})\bigr)^r\;P_{:t}(\mathrm{d}\xi_{:t})\right)^{\nicefrac1r}.
	\end{equation}
	The above estimate gives only an upper bound and equality does not hold in general, even when the optimal transport maps are inserted. To find appropriate transport mappings in dimensions $m > 1$ is not an easy job.  However, this estimate does not require the validity of a property like~\eqref{eq:cond}.
\end{proposition}

\begin{remark}
	For scalar processes ($\xi_t\in \mathbb R^1$, i.e., $m_t=1$) with $r=1$, the optimal transport map is the \emph{quantile coupling}: $S_{t,\xi_{:t-1}} = \tilde F_t^{-1}(\,\cdot\mid\xi_{:t-1})\circ F_t(\,\cdot\mid\xi_{:t-1})$, where~$F_t$ and~$\tilde F_t$ are the conditional distribution functions of~$P_t$ and~$\tilde P_t$ respectively.
	A variant of Vallander's formula then gives
	\[	w_r\bigl(P_t(\,\cdot\mid\xi_{:t-1}),\,\tilde P_t(\,\cdot\mid\xi_{:t-1})\bigr)= \Bigl(\int_0^1 \bigl|F_t^{-1}(u\mid\xi_{:t-1}) - \tilde F_t^{-1}(u\mid\xi_{:t-1})\bigr|^r\,\mathrm{d}u\Bigr)^{1/r}, \]
	making~\eqref{eq:6} fully explicit.
\end{remark}

\begin{algorithm}[t]
	\caption{Generation of a quantized scenario tree}\label{alg:QuantTree}
	\KwIn{Law~$P$ of a stochastic process; accuracy~$\delta>0$;
		branching numbers~$n_1,\dots,n_T$ from Proposition~\ref{prop:branches} below}
	\BlankLine
	Compute an optimal $n_1$-quantization of~$P_1$: find support points $\tilde\xi_1^1,\dots,\tilde\xi_1^{n_1}\in\Xi_1$ and weights
	\[	\tilde p_1^j = P_1(V_1^j),\] where $V_1^j$ is the Voronoi cell of $\tilde\xi_1^j$ under~$w_r$\;
	\For{$t = 2, \dots, T$}{
		\ForEach{node $\tilde\xi_{:t-1}$ at stage $t-1$}{
			Compute an optimal $n_t$\nobreakdash-quantization of the conditional kernel~$P_t(\,\cdot\mid\tilde\xi_{:t-1})$: find support points $\tilde\xi_t^1,\dots,\tilde\xi_t^{n_t}\in\Xi_t$ and weights \[	\tilde p_t^j = P_t(V_t^j\mid\tilde\xi_{:t-1}),\]
			where $V_t^j$ is the corresponding Voronoi cell.
	}	}
	\KwOut{Quantized scenario tree~$\tilde P = \tilde P_1(\mathrm{d}\tilde\xi_1)\cdots \tilde P_T(\mathrm{d}\tilde\xi_T\mid\tilde\xi_{:T-1})$, where
		\[	\tilde P_t(\,\cdot\mid\tilde\xi_{:t-1}) = \sum_{j=1}^{n_t}\tilde p_t^j\,\delta_{\tilde\xi_t^j}.\]}
\end{algorithm}
Proposition~\ref{prop:Delta} reduces verification of condition~\eqref{eq:cond} to a Lipschitz property of the true process and a uniform approximation error of the conditional kernels.
For this, approximations for $\varepsilon_t$ have to be found which can be established with Zador’s theorem.

With Zador’s theorem (cf.~\eqref{eq:Zador} above), we have that
\[	e_r(n_t)\coloneqq \inf_{|\supp(\tilde{P})|= n_t}w_r\bigl( P_t(\,\cdot\mid|\xi_{:t-1}),\ \tilde{P}\bigr)\sim q_t\cdot n_t^{-\nicefrac{1}{m_t}}, \]
where $q_t\coloneqq Q_{m_t,r}\cdot  \sup_{u_{:t-1}} \|f_t(\,\cdot\mid u_{:t-1})\|^{\nicefrac1r}_{\nicefrac{m_t}{(m_t+r)}}$ is Zador’s (dimension--dependent) constant of the kernel at~$u_{:t-1}$. Here $f_t(\,\cdot\mid u_{:t-1})$ is the conditional density of $P_t(\,\cdot\mid u_{:t-1})$ given $u_{:t-1}$.
With~\eqref{eq:5}, it follows that 
\[	\nd_r(P,\tilde{P})\le \sum_{t=1}^T q_t\cdot n_t^{-\nicefrac{1}{m_t}}\prod_{j=t+1}^T\kappa_{j}. \]
If the process does not expand over time (i.e., $\kappa_t= 1$), then the nested distance is bounded by 
\[	\nd_r(P,\tilde{P})\le \sum^T_{t=1}q_t\cdot n_t^{-\nicefrac{1}{m_t}}. \]
So to guarantee that $\nd_r(P,\tilde{P})\le\delta,$ it is enough to choose a branching structure $(n_{1},n_{2},\dots,n_T)$ so that
\begin{align}
	\sum_{t=1}^T q_t\cdot n_t^{-\nicefrac{1}{m_t}} \prod_{j=t+1}^T\kappa_j\le \delta.\label{eq:b1}
\end{align}
We now determine the optimal branching numbers $(n_1,\dots,n_T)$ that minimize the total scenario count subject to the accuracy constraint~\eqref{eq:b1}.

\subsection{Minimizing the number of leaves}
The total number of scenarios, which correspond to the leaves in the quantization tree with~$n_t$ branchings at each node at stage~$t$, is 
\begin{equation}
	N= \prod_{t=1}^T n_t.\label{eq:Leaves}
\end{equation}
We shall refer to elementary stochastic processes with $n_t$ branches at each node at stage~$t$ by their branching structure $(n_{1},n_{2},\dots,n_T)$. 

In what follows, we are interested in achieving an approximation of a measure $P$ by an elementary tree~$\tilde{P}$ so that the approximation in nested distance is $\nd_r(P,\tilde{P})\le \delta$ ($\delta>0$) while minimizing the total number of scenarios~$N$. Notice that we work with an upper bound of $\nd_r$, the true value may be smaller.
\begin{proposition}\label{prop:branches}
	Asymptotically, an approximation error of~$\delta$ in nested distance can be achieved with a tree with
	\[
		N= \prod_{t=1}^T\left(\frac{m\,q_t}{m_t\delta}K_t\right)^{m_t}
	\]
	leaves, where $m=\sum^T_{t=1}m_t$ is the total dimension and $K_t\coloneqq \prod_{t^\prime>t}\kappa_{t^\prime}$.
	The branchings at stage~$t$ are given explicitly by 
	\begin{equation}
		n_t= \left\lceil \left(\frac{m\,q_t K_t}{m_t}\right)^{m_t}\delta^{-m_t} \right\rceil\label{eq:111}
	\end{equation}
	provided that the root of the tree equals $\xi_0$, the starting value of the process $\xi$. Here, $\lceil z \rceil$ is the integer ceiling of $z$. 
\end{proposition}
\begin{proof}
	We first derive the optimal branching structure.
	To minimize $\log N=\sum^T_{t=1}\log n_{t}$ subject to the constraint $\sum_{t=1}^T q_t K_t n_t^{-\nicefrac{1}{m_t}}= \delta$ (cf.~\eqref{eq:5}), the Lagrangian function with Lagrange multiplier~$\lambda$ is
	\[ L= \sum_{t=1}^T\log n_t+ \lambda\left(\sum_{t=1}^T q_t K_t n_t^{-\nicefrac{1}{m_t}}-\delta\right). \]
	Differentiating with respect to $n_t$ gives the Lagrangian stationary conditions
	\[
		\frac{1}{n_t}-\lambda\frac{1}{m_t}q_t K_t n_t^{-\frac{1}{m_t}-1}=0
	\]
	so that 
	\begin{equation}
		n_t=\left(\lambda\frac{q_t K_t}{m_t}\right)^{m_t}.\label{eq:112}
	\end{equation}
	Imposing the budget constraint fixes $\lambda$, giving 
	\[
		\delta= \sum_{t=1}^T q_t K_t\left(\frac{m_t}{\lambda\,q_t K_t}\right)= \frac{1}{\lambda}\sum_{t=1}^T m_t,
	\]
	so that $\lambda=\frac{m}{\delta}$ with $m= \sum_{t=1}^T m_t$.
	With this, the branching~\eqref{eq:112} is $n_t= \left(\frac{m\,q_t K_t}{m_t}\right)^{m_t}\delta^{-m_t}$.
	It follows that the branching structure~\eqref{eq:111} is optimal.

	Note as well that 
	\[	\prod_{t=1}^T n_t= N \]
	and with~\eqref{eq:5} and \(\varepsilon_t \coloneqq q_t n_t^{-\nicefrac{1}{m_t}}\), it holds that
	\begin{align*}
		\nd_r(P,\tilde{P})& \le \sum_{t=1}^T q_t\cdot n^{-\nicefrac{1}{m_t}}_t\cdot\prod_{j=t+1}^T\kappa_{j}\\
			& =\sum_{t=1}^T q_t\cdot\frac{m_t}{m\,q_t\,\prod_{j=t+1}^T\kappa_{j}}\delta\cdot\prod^T_{j=t+1}\kappa_{j}\\
			& =\sum_{t=1}^T\frac{m_t}{m}\delta=\delta,
	\end{align*}
	which is the desired assertion.
\end{proof}

\section{Random trees}\label{sec:random}
The random scenario tree is built stage by stage.
At each stage~$t$ and each node reached, the conditional distribution of~$\xi_t$ given the path history is approximated by drawing~$n_t$ independent samples from that conditional.
We shall denote the measure governing the independent samples of~$P$ by
\[	\mathds P\coloneqq P\otimes P\otimes\dots.\]

Each draw at a node $\tilde\xi_{:t-1}$ implicitly defines a stage-wise transport map $S_{t,\tilde\xi_{:t-1}}$ pushing $P_t(\,\cdot\mid\tilde\xi_{:t-1})$ to the empirical measure $\tilde P_t(\,\cdot\mid\tilde\xi_{:t-1})$.
Composing these maps consistently (cf.\ Proposition~\ref{prop:mapping}) yields the mapping-based viewpoint of the random tree: the nested distance error is the $P$-integrated transport cost $\bigl(\int_\Xi d(\xi,S_T(\xi))^r\,P(\mathrm{d}\xi)\bigr)^{\nicefrac1r}$, decomposed stage by stage via~\eqref{eq:6}. 
Algorithm~\ref{alg:Tree} summarizes the concrete steps to construct the random tree.

Approximating the tree at each node by samples as in~\eqref{eq:8} causes an approximation error.
At each node $\tilde\xi_{:t-1}$, the approximation error for the $t$-th stage kernel is the Wasserstein distance
\begin{equation}\label{eq:7}
	\mathcal{E}_t(\tilde\xi_{:t-1})\coloneqq w_r\bigl(P_t(\,\cdot\mid\tilde\xi_{:t-1}),\;\tilde P_t(\,\cdot\mid\tilde\xi_{:t-1})\bigr),
\end{equation}
which is a random variable depending on the~$n_t$ samples drawn at that node.
By Proposition~\ref{prop:Delta} and Theorem~\ref{thm:Mirkov}, the overall nested distance is controlled by the worst-case stage errors $\varepsilon_t= \sup_{\tilde\xi_{:t-1}}\mathcal{E}_t(\tilde\xi_{:t-1})$.
\begin{algorithm}[t]
	\caption{Generation of a random tree}\label{alg:Tree}
	\KwIn{Law~$P$ of a stochastic process; branching numbers $n_1,\dots,n_T$}
	\KwOut{Random tree~$\tilde P$ with $N=\prod_{t=1}^T n_t$ leaves}
	\BlankLine
	Draw $n_1$ independent realizations $\tilde\xi_1^{(1)},\dots,\tilde\xi_1^{(n_1)}$ from the marginal~$P_1$; these form the $N_1=n_1$ nodes at stage~$1$ with empirical measure
	\[	\tilde P_1 = \frac{1}{n_1}\sum_{i_1=1}^{n_1}\delta_{\tilde\xi_1^{(i_1)}}. \]
	\For{$t = 2,\dots,T$}{
		\ForEach{node $\tilde\xi_{:t-1}=\bigl(\tilde\xi_1^{(i_1)},\dots,\tilde\xi_{t-1}^{(i_{t-1})}\bigr)$
			among the $N_{t-1}=\prod_{s<t}n_s$ nodes at stage $t-1$}{
			Draw $n_t$ independent realizations
			\begin{equation}\label{eq:8}
				\tilde\xi_t^{(i_1,\dots,i_{t-1},1)},\;\dots,\;\tilde\xi_t^{(i_1,\dots,i_{t-1},n_t)}
			\end{equation}
			from the kernel $P_t(\,\cdot\mid\tilde\xi_{:t-1})$\;
			Set these $n_t$ realizations as the children of $\tilde\xi_{:t-1}$ with empirical conditional measure
			\[	\tilde P_t\bigl(\,\cdot\mid\tilde\xi_{:t-1}\bigr)= \frac{1}{n_t}\sum_{j=1}^{n_t}\delta_{\tilde\xi_t^{(i_1,\dots,i_{t-1},j)}}\;. \]
		}	}
	\BlankLine
	\KwOut{The tree measure $\tilde P$ has disintegration (cf.~\eqref{eq:disintegrated})
	\[	\tilde P(\mathrm{d}\tilde\xi_1,\dots,\mathrm{d}\tilde\xi_T)= \tilde P_1(\mathrm{d}\tilde\xi_1)\, \tilde P_2(\mathrm{d}\tilde\xi_2\mid\tilde\xi_1)\,\cdots\, \tilde P_T(\mathrm{d}\tilde\xi_T\mid\tilde\xi_{:T-1}). \]
	Each path $\bigl(\tilde\xi_1^{(i_1)},\dots,\tilde\xi_T^{(i_1,\dots,i_T)}\bigr)$ carries equal probability $1/N$}
\end{algorithm}

Probabilistic bounds for stochastic programming approximations via large-deviation inequalities were derived by \citet{Wets1995}; the present approach uses the sharper sub-Gaussian concentration of \citet{Bolley}, for which we recall the following theorem.
\begin{theorem}\label{thm:Bolley}
	Let~$P$ be a probability measure on $\Xi= \mathbb{R}^{m}$ with $\int_\Xi e^{\alpha\|u\|^{2}}P(\mathrm{d}u)< \infty$ for some $\alpha>0$, and let $r\in[1,2)$. Let $\tilde{P}_n=\frac{1}{n}\sum_{i=1}^n\delta_{X_{i}}$ with $X_{i}$ i.i.d.\ with respect to~$P$.
	Then for every $m^\prime>m$ and $c<\alpha$ there is $n_{0}$ such that for all $\varepsilon>0$ and
	$n\geq n_{0}\max(1,\ \varepsilon^{-(m^\prime+2)})$,
	\[	P\bigl(w_r(P,\tilde{P}_n)>\varepsilon\bigr)\leq e^{-c\,n\varepsilon^{2}}. \]
\end{theorem}

The bound in Theorem~\ref{thm:Bolley} holds for a single distribution. In the multistage setting, Algorithm~\ref{alg:Tree} draws $n_t$ independent samples from \emph{each} conditional distribution $P_t(\,\cdot\mid\tilde\xi_{:t-1})$ at every node $\tilde\xi_{:t-1}$ of stage~$t-1$.
We aggregate the resulting concentration bounds over all nodes and all stages by a union bound, then invoke Theorem~\ref{thm:Mirkov} to pass from conditional Wasserstein errors to the nested distance.

We assume throughout that all conditional distributions $P_t(\,\cdot\mid\xi_{:t-1})$, for $t=1,\ldots,T$ and all histories $\xi_{:t-1}$, satisfy $\sup_{\xi_{:t-1}} \int_{\Xi_t}e^{\alpha\|u\|^2}P_t(\mathrm{d}u\mid\xi_{:t-1})< \infty$ for some 
$\alpha>0$.

Let $c_t>0$ be the concentration constant from Theorem~\ref{thm:Bolley} for stage~$t$, and set $c\coloneqq \min_t c_t$. As before, we assume that the root of the tree~$\tilde{\xi}_0$ equals the starting value $\xi_0$ of the scenario process.

Suppose the true process $P$ satisfies the Lipschitz condition
\begin{equation}\label{eq:Lip}
	w_r\bigl(P_{t+1}(\,\cdot\mid\xi_{:t}),\;P_{t+1}(\,\cdot\mid\tilde{\xi}_{:t})\bigr) \le \kappa_t\cdot D_t(\xi_{:t},\tilde{\xi}_{:t})
\end{equation}
for all $\xi_{:t},\tilde{\xi}_{:t}\in\Xi_{:t}$ and each stage $t$.
By Proposition~\ref{prop:Delta}, this implies condition~\eqref{eq:cond} at every stage, with
\[	\varepsilon_{t+1} \coloneqq \sup_{\tilde{\xi}_{:t}} w_r\bigl(P_{t+1}(\,\cdot\mid\tilde{\xi}_{:t}),\; \tilde{P}_{t+1}(\,\cdot\mid\tilde{\xi}_{:t})\bigr).\]

\begin{theorem}[Nested distance of the random tree]\label{thm:random}
	Let $e_1,\ldots,e_T>0$, $\eta\in(0,1)$, and $r\in[1,2)$.
	Suppose the branching numbers $n_1,\ldots,n_T$ in Algorithm~\ref{alg:Tree} satisfy, for each $t=1,\ldots,T$,
	\begin{equation}\label{eq:ntcond}
		n_t \ge \frac{\log N_{t-1}+\log(T/\eta)}{c\,e_t^2} \qquad\text{and}\qquad
		n_t \ge n_0\,\max\!\bigl(1,\,e_t^{-(m_t^\prime+2)}\bigr),
	\end{equation}
	where $N_{t-1}=\prod_{j=1}^{t-1}n_j$ is the number of nodes at stage~$t-1$ (with $N_0=1$) and $m_t^\prime>m_t$.  Then
	\begin{equation}\label{eq:mainbound}
		\mathds P\left(\nd_r(P,\tilde P) > \sum_{t=1}^T e_t\,K_t\right) \le\eta
	\end{equation}
with $K_t\coloneqq \prod_{s > t} \kappa_s$.
\end{theorem}

\begin{proof}
	Fix a stage $t$ and a node $\tilde\xi_{:t-1}$ at stage~$t-1$.
	Algorithm~\ref{alg:Tree} produces $n_t$ i.i.d.\ samples from $P_t(\,\cdot\mid\tilde\xi_{:t-1})$, yielding the empirical measure $\tilde P_t(\,\cdot\mid\tilde\xi_{:t-1})$.
	By Theorem~\ref{thm:Bolley} and the second condition in~\eqref{eq:ntcond},
	\begin{equation*}
		\mathds P\left(w_r\bigl(P_t(\,\cdot\mid\tilde\xi_{:t-1}),\, \tilde P_t(\,\cdot\mid\tilde\xi_{:t-1})\bigr)>e_t\right) \leq e^{-c\,n_t\,e_t^2}.
	\end{equation*}
	The first condition in~\eqref{eq:ntcond} ensures $e^{-c\,n_t e_t^2}\le \frac{\eta}{N_{t-1}T}$.
	A union bound over all $N_{t-1}$ nodes at stage~$t-1$ gives
	\begin{equation*}
		\mathds P\left(\varepsilon_t>e_t\right) \le N_{t-1}\cdot\frac{\eta}{N_{t-1}\,T} =\frac{\eta}{T},
	\end{equation*}
	where $\varepsilon_t= \sup_{\tilde\xi_{:t-1}} w_r(P_t(\,\cdot\mid\tilde\xi_{:t-1}),\tilde P_t(\,\cdot\mid\tilde\xi_{:t-1}))$ is the maximal conditional Wasserstein error at stage~$t$ (cf.\ Proposition~\ref{prop:Delta}).
	A further union bound over $t=1,\ldots,T$ yields
	\begin{equation*}
		\mathds P\left(\exists\,t\colon \varepsilon_t>e_t\right)\le\eta.
	\end{equation*}
	On the complementary event $\{\varepsilon_t\leq e_t\text{ for all }t\}$, Proposition~\ref{prop:Delta} verifies condition~\eqref{eq:cond} with parameters $\varepsilon_t\le e_t$, so Theorem~\ref{thm:Mirkov} gives $\nd_r(P,\tilde P)\le \sum_{t=1}^T \varepsilon_t\,K_t\le \sum_{t=1}^T e_t\,K_t$, which is~\eqref{eq:mainbound}.
\end{proof}

It remains to choose the error budget $(e_t)$ so that $\sum_{t=1}^T e_t K_t= \delta$ and the total number of scenarios $N= \prod_{t=1}^T n_t$ is as small as possible.

\subsection{Acceptable error budget}
The following proposition allocates the errors to not exceed the total budget of~$\delta$.

\begin{proposition}[Optimal error allocation]\label{prop:erroralloc}
	Fix $\delta>0$.  Among all allocations satisfying $\sum_{t=1}^T e_t K_t= \delta$, the one minimizing $\log N=\sum_t\log n_t$ (up to logarithmic factors in $\delta^{-1}$ and $\eta^{-1}$) is the equal-contribution rule
	\begin{equation}\label{eq:optalloc}
		e_t= \frac{\delta}{T\,K_t}, \qquad t=1,\ldots,T.
	\end{equation}
	With this choice the conditions~\eqref{eq:ntcond} are met by
	\begin{equation}\label{eq:optbranch}
		n_t= \left\lceil\frac{T^2 K_t^2}{c\,\delta^2} \bigl(\log N_{t-1}+\log(T/\eta)\bigr)\right\rceil,
	\end{equation}
	and the total number of scenarios satisfies $N_\delta= O\left(\delta^{-2T} \bigl(\log(\delta^{-1})\bigr)^T\right)$.
\end{proposition}

\begin{proof}
	Ignoring logarithmic factors, the first condition in~\eqref{eq:ntcond} gives $n_t\asymp e_t^{-2}$, so $\log N\asymp -2\sum_t\log e_t$.
	Minimizing this over $(e_t)$ subject to $\sum_t e_t K_t=\delta$ by Lagrange multipliers yields $e_t\propto K_t^{-1}$; the constraint then fixes the constant, giving \eqref{eq:optalloc}.
	Substituting into the first bound of~\eqref{eq:ntcond} produces~\eqref{eq:optbranch}.
	With this choice $n_t=\mathcal O(K_t^2\delta^{-2}\log(\delta^{-1}))$, so $N=\prod_t n_t= \mathcal{O}\bigl(\delta^{-2T}(\log\delta^{-1})^T\bigr)$, the assertion.
\end{proof}
\begin{corollary}[Probabilistic nested-distance guarantee]\label{cor:main}
	Under the hypotheses of Theorem~\ref{thm:random}, choose the equal-contribution allocation $e_t = \delta/(TK_t)$ and the branching numbers~\eqref{eq:optbranch}. Then the random tree produced by Algorithm~\ref{alg:Tree} satisfies
	\begin{equation}\label{eq:mainresult}
		\mathds P\left(\nd_r(P,\tilde{P}) > \delta\right) \leq \eta,
	\end{equation}
	and the total number of leaves satisfies $N_\delta = \mathcal O\left(\delta^{-2T}(\log\delta^{-1})^T\right)$.
\end{corollary}
\begin{proof}
  The allocation satisfies $\sum_{t=1}^T e_t K_t = \delta$, so Theorem~\ref{thm:random} gives~\eqref{eq:mainresult}.
  The scenario count is Proposition~\ref{prop:erroralloc}.
\end{proof}

\begin{remark}[Comparison with quantized trees]\label{rem:compare}
	Propositions~\ref{prop:branches} and~\ref{prop:erroralloc} both achieve nested-distance accuracy~$\delta$, but with different scenario counts and guarantees.

	\emph{Quantized tree (Proposition~\ref{prop:branches}).}  Zador's theorem gives $n_t\asymp e_t^{-m_t}$, which after the same Lagrangian optimization yields $N_{\mathrm{quant}}= \mathcal O(\delta^{-m})$ with $m=\sum_{t=1}^Tm_t$.
	The bound $\nd_r(P,\tilde P)\le\delta$ holds \emph{deterministically}.
	
	\emph{Random tree (Proposition~\ref{prop:erroralloc}).}  Theorem~\ref{thm:Bolley} gives $n_t\asymp e_t^{-2}$ (independently of dimension), so $N_{\mathrm{rand}}= \mathcal O(\delta^{-2T}\,\mathrm{polylog})$.
	The bound $\nd_r(P,\tilde P)\le\delta$ holds with probability at least $1-\eta$.
	
	In the scalar case $m_t=1$ one has $m=T$, so
	$N_{\mathrm{quant}}= \mathcal O(\delta^{-T})$ against $N_{\mathrm{rand}}= \mathcal O(\delta^{-2T})$:
		random sampling requires quadratically more scenarios for a weaker, probabilistic guarantee.  The picture reverses for high-dimensional stages: whenever $m_t>2$, the Bolley rate $e_t^{-2}$ is faster than the quantization rate $e_t^{-m_t}$, so the random tree is more parsimonious at those stages.
		Additionally, generating the random tree requires only i.i.d.\ sampling, whereas computing an optimal quantization is an NP-hard problem in general dimension.
\end{remark}

\subsection{Confidence bound}
Accepting an error of $\delta$, the confidence in~\eqref{eq:mainresult} improves.
The following corollary elaborates the dependency explicitly.
\begin{corollary}[Dependence on the confidence level~$\eta$]\label{cor:eta}
	Under the hypotheses and notation of Proposition~\ref{prop:erroralloc}, fix $\delta>0$.  As $\eta\to 0$, the equal-contribution allocation requires
	\[	N = \mathcal O\left(\Bigl(\log\tfrac{1}{\eta}\Bigr)^T\right) \]
 	leaves.  In particular, the dependence on the confidence level~$\eta$ is merely polylogarithmic.
	If~$\eta$ and~$\delta$ tend to zero simultaneously with $\eta \geq \delta^\alpha$ for some fixed $\alpha>0$, then $\log(1/\eta)\leq \alpha\log(1/\delta)$ and the bound $N= \mathcal O\left(\delta^{-2T}(\log\delta^{-1})^T\right)$ of Proposition~\ref{prop:erroralloc} continues to hold.
\end{corollary}
\begin{proof}
	Set $\ell\coloneqq\log(1/\eta)$. With $e_t= \delta/(TK_t)$, the branching numbers in Theorem~\ref{thm:random} satisfy
	\[	n_t= \Bigl\lceil \tfrac{T^2K_t^2}{c\,\delta^2} \bigl(\log N_{t-1}+\log(T/\eta)\bigr) \Bigr\rceil. \]
	We show by induction on $t$ that $N_{t-1}= \mathcal O(\ell^{t-1})$.

	For $t=1$ this holds with $N_0=1$.
	Assuming $N_{t-1}= \mathcal O(\ell^{t-1})$, we get $\log N_{t-1}= \mathcal O(\log\ell)=o(\ell)$, so $n_t= \mathcal O(\ell)$, and $N_t=N_{t-1}\cdot n_t= \mathcal O(\ell^t)$.
	At the final stage $t=T$ this gives $N=\mathcal O(\ell^T)$. The second claim follows because $\log(1/\eta)\leq\alpha\log(1/\delta)$ whenever $\eta\geq\delta^\alpha$.
\end{proof}
Corollary~\ref{cor:eta} reveals an important asymmetry between the two approximation parameters~$\delta$ and~$\eta$.
Achieving higher accuracy --- that is, reducing $\delta$ --- is expensive: the tree size grows as $\delta^{-2T}$, a power law whose exponent scales with the number of stages.
Achieving higher confidence --- that is, reducing $\eta$ --- is essentially free by comparison: the tree size grows only as $(\log(1/\eta))^T$, so even small failure probabilities require only a moderate increase in~$N$.
For instance, with $T=3$ stages, reducing $\eta$ from $5\%$ to $10^{-10}$ multiplies $\log(1/\eta)$ by a factor of roughly $10/3$, and hence $N$ by $(10/3)^3\approx 37$ --- a modest constant, independent of $\delta$.
By contrast, halving~$\delta$ multiplies~$N$ by $2^{2T}=64$.

This behavior is characteristic of Monte Carlo methods, and it reflects the origin of the $\log(1/\eta)$ term: it enters through a union bound over the nodes of the tree, which contributes $\log(T/\eta)$ to each branching number --- a cost that is linear in $\log(1/\eta)$ but independent of the accuracy $\delta$. The accuracy $\delta$, on the other hand, enters through the quantization error $e_t\propto\delta$, which appears as $e_t^{-2}$ in the sample-size condition of the Bolley-Guillin-Villani theorem, driving the polynomial rate.

In practical terms, the random tree construction is therefore well suited to problems where high confidence is required but the approximation accuracy~$\delta$ is not too small.
When $\delta$ is very small, the $\delta^{-2T}$~factor dominates, and the deterministic quantized tree, whose size is~$O(\delta^{-m})$, becomes preferable whenever the total dimension $m= \sum_t m_t$ satisfies $m<2T$.

\section{Conclusion}
We derive quantitative bounds on the nested distance between a stochastic process and two families of finite scenario-tree approximations, and identify the optimal branching structure for each family.

For the \emph{deterministic quantization tree} (Algorithm~\ref{alg:QuantTree}), Proposition~\ref{prop:branches} shows that the nested distance bound $\nd_r(P,\tilde P)\le\delta$ can be guaranteed with certainty using
\[
    N= \prod_{t=1}^T\left(\frac{m\,q_t K_t}{m_t\,\delta}\right)^{m_t} = \mathcal{O}(\delta^{-m})
\]
leaves, where $m= \sum_t m_t$ is the total state-space dimension and $K_t= \prod_{j>t}\kappa_j$ captures the expansion of the process. The optimal branching numbers follow explicitly from Zador's quantization coefficients and the expansion parameters.

For the \emph{random sampling tree} (Algorithm~\ref{alg:Tree}),
Theorem~\ref{thm:random} shows that $\nd_r(P,\tilde P)\le\delta$ holds with probability at least $1-\eta$ using
\[
    N= \mathcal{O}\bigl(\delta^{-2T}(\log\delta^{-1})^T\bigr)
\]
leaves.
This count is independent of the ambient dimensions $m_t$; the exponent in the 
scenario count is governed by the number of stages~$T$ alone, not by the 
state-space dimensions.
The sub-Gaussian concentration condition for processes with sub-Gaussian tails is the key structural assumption.

The two constructions offer complementary guarantees. The deterministic tree provides an exact bound, but its cost scales with the total dimension $m=\sum_t m_t$, making it prohibitive in high-dimensional settings. The random tree provides a probabilistic bound at a dimension-free cost: for processes with large state-space dimension ($m>2T$) and an acceptable failure probability~$\eta$, it achieves the target nested distance bound with substantially fewer leaves.

Extensions to non-product spaces, to heavy-tailed conditional distributions, and to state-dependent expansion coefficients are natural. From an algorithmic perspective, the optimal quantization step in Algorithm~\ref{alg:QuantTree} relies on efficient vector-quantization solvers whose computational cost depends on both the dimension~$m_t$ and the kernel, and warrants further study.
Finally, translating the nested-distance bounds derived here into explicit error bounds on the optimal value and optimal policy of the multistage problem~\eqref{eq:1} would connect the approximation-theoretic results of this paper directly to algorithmic practice, cf.\ \citep{DupachovaWets1988, Romisch2003}.

\section*{Acknowledgment}
The authors used Claude (Anthropic) for assistance with manuscript preparation, writing, and editing. All mathematical content, results, and proofs are solely the work of the authors.

\label{References}
\bibliographystyle{abbrvnat}
\bibliography{LiteraturAlois}

\end{document}